\documentstyle[amscd,amssymb,verbatim,12pt]{amsart}
\newcommand{\ch}{\mbox{\rm ch}}

\newcommand{\Ker}{\mbox{\rm Ker}}

\newcommand{\Image}{\mbox{\rm Im}}
\newcommand{\Aut}{\mbox{\rm Aut}}
\newcommand{\tr}{\mbox{\rm tr}}
\newcommand{\Tr}{\mbox{\rm Tr}}
\newcommand{\TR}{\mbox{\rm TR}}

\newcommand{\End}{\mbox{\rm End}}

\newcommand{\HH}{\mbox{\rm H}}
\newcommand{\HC}{\mbox{\rm HC}}

\newcommand{\vol}{\mbox{\rm vol}}
\newcommand{\Hom}{\mbox{\rm Hom}}

\newcommand{\sign}{\mbox{\rm sign}}

\newcommand{\spec}{\mbox{\rm spec}}

\newcommand{\Diff}{\mbox{\rm Diff}}

\newcommand{\R}{{\Bbb R}}
\newcommand{\C}{{\Bbb C}}

\newcommand{\Z}{{\Bbb Z}}
\newcommand{\htriangle}{\widehat{\triangle}}
\newcommand{\hD}{\widehat{D}}
\newcommand{\hGamma}{\widehat{\Gamma}}
\newcommand{\hM}{\widehat{M}}
\newcommand{\hE}{\widehat{E}}
\newcommand{\hZ}{\widehat{Z}}
\newcommand{\hF}{\widehat{F}}
\newcommand{\hphi}{\widehat{\phi}}
\theoremstyle{plain}
\newtheorem{definition}{Definition}

\newtheorem{lemma}{Lemma}

\newtheorem{proposition}{Proposition}

\newtheorem{conjecture}{Conjecture}

\numberwithin{equation}{section}
\renewcommand{\rm}{\normalshape}
\begin{document}
\title{Delocalized $L^2$-Invariants}
\author{John Lott}
\address{Department of Mathematics\\
University of Michigan\\
Ann Arbor, MI  48109-1109\\
USA}
\email{lott@@math.lsa.umich.edu}
\thanks{Research supported by NSF grant DMS-9403652}
\date{November 24, 1996}
\maketitle
\begin{abstract}
We define extensions of the $L^2$-analytic invariants of closed manifolds,
called delocalized $L^2$-invariants.
These delocalized invariants are constructed in terms of 
a nontrivial conjugacy class of the fundamental 
group. We show that in many cases, they are topological in nature.  We show
that the marked length spectrum of an odd-dimensional hyperbolic manifold
can be recovered from its delocalized $L^2$-analytic torsion. There are
technical convergence questions.
\end{abstract}
\section{Introduction}
Let $M^d$ be a closed $d$-dimensional
hyperbolic manifold. Within the class of such
manifolds, the volume $\vol(M)$ of $M$ is a topological invariant; 
this follows from Mostow rigidity \cite{Mostow (1973)} if $d > 2$ and from the
Gauss-Bonnet theorem if $d = 2$.

One can ask whether there is a topological
invariant defined on a wider class of smooth
manifolds which, in the case of a hyperbolic manifold, reproduces its
hyperbolic volume.  
One such invariant is Gromov's simplicial volume
\cite{Gromov (1982)}, which is defined for all compact manifolds $M$.  

If $d$ is even then another relevant topological invariant is the Euler 
characteristic $\chi(M)$. 
If $M$ is hyperbolic then $\vol(M) = (-1)^{d/2} \: \frac{vol(S^d)}{2} \: 
\chi(M)$. Equivalently, we could phrase this in terms of
the $L^2$-Euler characteristic $\chi_{(2)}(M)$, which equals $\chi(M)$
\cite{Atiyah (1976)}.

An odd-dimensional counterpart is
the $L^2$-analytic torsion $T_{\langle e \rangle}(M)$ 
\cite{Lott (1992b),Mathai (1992)}, which is a (smooth) topological invariant of
compact manifolds $M$ with vanishing $L^2$-Betti numbers and positive
Novikov-Shubin invariants. 
Odd-dimensional closed hyperbolic manifolds satisfy
these conditions.  For such manifolds,
$T_{{\langle e \rangle}}(M) = c_d \: \vol(M)$
where $c_d$ only depends on $d$.
For example, $c_3 = - \frac{1}{3\pi}$.

Besides the volume,
another interesting invariant of a closed hyperbolic manifold is its marked
length spectrum.  This is the function which, to a conjugacy class
$\langle g \rangle \subset \pi_1(M)$, assigns the hyperbolic length
of the unique closed geodesic in the free
homotopy class specified by $\langle g \rangle$. In this paper we address
the question of whether there is a topological invariant of some class of
smooth manifolds which, in the case of a hyperbolic manifold of dimension
greater than two, reproduces its
marked length spectrum.

We consider three invariants of a closed Riemannian manifold $M$, 
called the delocalized $L^2$-Betti numbers, 
$L^2$-analytic torsion and $L^2$-eta invariant. 
These invariants are defined in Section 
\ref{Definitions and Statements of Results}.
The reason for the word ``delocalized'' is that the ordinary
$L^2$-invariants 
\cite{Atiyah (1976),Cheeger-Gromov (1985),Lott (1992b),Mathai (1992)}
are localized at the
identity element of the fundamental group, in a sense which will be made
precise. 
In contrast, the delocalized invariants are defined in terms of a nontrivial
conjugacy class $\langle g \rangle$ of the fundamental group.

The delocalized $L^2$-invariants are constructed using
the heat kernel on the universal cover $\widetilde{M}$ of $M$. 
There are technical difficulties in showing that the formal expressions
for the invariants are actually well-defined. These difficulties involve
estimating heat kernels on $\widetilde{M}$ at large-distance and
large-time simulaneously.
We cannot show that the formal expressions make sense in full generality.  
However, we do not know of any examples in which this is not the case.
We do show that the delocalized $L^2$-invariants are well-defined and
metric-independent in the following cases :\\
1. If $\langle g \rangle$ is a finite conjugacy class.\\
2. If the 
fundamental group is virtually nilpotent or Gromov-hyperbolic and if there is
a gap away from zero in the spectrum of the lifted $p$-form Laplacian, in the
case of the delocalized $L^2$-Betti number.\\
3. If the 
fundamental group is virtually nilpotent or Gromov-hyperbolic and if there is
a gap away from zero in the spectrum of the lifted Dirac-type operator, in the
case of the delocalized $L^2$-eta invariant.\\
 
We compute the delocalized $L^2$-invariants in some cases of interest.
For an odd-dimensional hyperbolic manifold, we show that the delocalized
$L^2$-invariants are well-defined and that the delocalized $L^2$-analytic 
torsion reproduces the marked length
spectrum. This is based on work of Fried
\cite{Fried (1986)} and Millson \cite{Millson (1978)}.
In the case of a mapping torus with simply-connected fiber,
we show that the delocalized $L^2$-invariants can be written in
terms of the action of the gluing diffeomorphism on the fiber cohomology.
In the case of a mapping torus 
whose fiber has finite fundamental
group, we show that the delocalized $L^2$-analytic torsion is determined by
the Nielsen fixed point indices of the gluing diffeomorphism.

If the fundamental group is torsion-free then the delocalized $L^2$-Betti
numbers vanish in all cases in which we can compute them.  We do not know
if this is always the case.

In order to really say that we have defined topological invariants
which reproduce the marked length spectrum of an odd-dimensional hyperbolic
manifold, it remains to prove the following.
\begin{conjecture}
Let $M$ be a closed odd-dimensional manifold which admits a hyperbolic
structure.  Then for any
Riemannian metric on $M$, the delocalized $L^2$-invariants are
well-defined and independent of the choice of metric.
\end{conjecture}

Another interesting question is whether there is a ``delocalized'' version
of the simplicial volume.

The structure of the paper is as follows.  In Section
\ref{Definitions and Statements of Results} we give the definitions of the
invariants and the
statements of the main results.  In Sections \ref{proof1}-\ref{proof8} 
we prove the results.  In Section \ref{Examples} we give examples to show 
that the results are not vacuous.

\section{Definitions and Statements of Results} 
\label{Definitions and Statements of Results}
Let $M^d$ be a closed connected oriented
Riemannian manifold. Let $\pi : \widehat{M} \rightarrow
M$ be a connected normal $\Gamma$-cover of $M$, equipped with
the pullback Riemannian metric. We let
$\gamma \in \Gamma$ act on $\widehat{M}$ on the right by $R_\gamma \in
\Diff(\widehat{M})$. Let ${\cal C}$ denote the set of conjugacy classes of 
$\Gamma$. 

\begin{definition} \label{conv}
Let ${\cal A}$ be the convolution algebra of elements 
$a \in C^\infty(\hM \times \hM)$ 
satisfying \\
1. $a(m \gamma, m^\prime \gamma) = a(m, m^\prime)$
for all $\gamma \in \Gamma$ and \\
2. There exists an $R_a > 0$ such that if $d(m, m^\prime) \ge R_a$ then
$a(m, m^\prime) = 0$.\\ \\
The multiplication in ${\cal A}$ is given by
\begin{equation}
(ab)(m, m^\prime) = \int_{\hM} a(m, m^{\prime \prime}) \:
b(m^{\prime \prime}, m^\prime) \: d\vol(m^{\prime \prime}).
\end{equation}

Given $a \in {\cal A}$ and $\langle g \rangle \in {\cal C}$, 
define $A_{\langle g \rangle} \in C^\infty(\hM)$ by
\begin{equation} \label{A} 
A_{\langle g \rangle}(m) = \sum_{\gamma \in \langle g \rangle} a(m
\gamma, m).
\end{equation}
\end{definition}

\begin{lemma}
For all $\gamma^\prime \in \Gamma$,
\begin{equation}
A_{\langle g \rangle}(m\gamma^\prime) = A_{\langle g \rangle}(m).
\end{equation}
\end{lemma}
\begin{pf}
We have
\begin{equation}
A_{\langle g \rangle}(m\gamma^\prime) = 
\sum_{\gamma \in \langle g \rangle} a(m\gamma^\prime \gamma, 
m \gamma^\prime) =
\sum_{\gamma \in \langle g \rangle} a \left(m\gamma^\prime \gamma 
\left(\gamma^\prime \right)^{-1}, m \right) = 
A_{\langle g \rangle}(m).
\end{equation}
\end{pf}

Thus $A_{\langle g \rangle} = \pi^* 
\alpha_{\langle g \rangle}$ for a unique $\alpha_{\langle g \rangle} \in
C^\infty(M)$.

\begin{definition}
Define $\Tr_{\langle g \rangle} : {\cal A} \rightarrow \C$ by
\begin{equation}
\Tr_{\langle g \rangle} (a) = \int_M \alpha_{\langle g \rangle} \: d\vol.
\end{equation}
\end{definition}

\begin{lemma} \label{comm}
For all $a, b \in {\cal A}$, 
\begin{equation}
\Tr_{\langle g \rangle} (ab) = \Tr_{\langle g \rangle} (ba).
\end{equation}
\end{lemma}
\begin{pf}
Let ${\cal F}$ be a fundamental domain for the
$\Gamma$-action on $\widehat{M}$. Formally,
\begin{align} \label{formally}
\Tr_{\langle g \rangle} (ab) & = \int_{\cal F} \sum_{\gamma \in \langle g 
\rangle} (ab)(m \gamma, m) \: d\vol(m) \\
& = \int_{\cal F} \sum_{\gamma \in \langle g 
\rangle} \int_{\hM} a(m\gamma, m^\prime) \: b(m^\prime, m) \: d\vol(m^\prime) 
\: d\vol(m) \notag \\
& = \int_{\cal F} \sum_{\gamma \in \langle g 
\rangle} \int_{\cal F} \sum_{\gamma^\prime \in \Gamma} 
a(m\gamma, m^\prime \gamma^\prime) \: b(m^\prime \gamma^\prime, m) \: 
d\vol(m^\prime) \: d\vol(m) \notag \\
& =  \sum_{\gamma \in \langle g 
\rangle} \sum_{\gamma^\prime \in \Gamma} \int_{\cal F} \int_{\cal F}  
b(m^\prime \gamma^\prime, m) \: a(m\gamma, m^\prime \gamma^\prime)  
\: d\vol(m^\prime)\: d\vol(m) \notag \\
& =  \sum_{\gamma \in \langle g 
\rangle} \sum_{\gamma^\prime \in \Gamma} \int_{\cal F} \int_{\cal F}  
b(m^\prime \gamma^\prime \gamma (\gamma^\prime)^{-1}, 
m\gamma (\gamma^\prime)^{-1}) \:
a(m\gamma (\gamma^\prime)^{-1}, m^\prime)
\: d\vol(m^\prime)\: d\vol(m) \notag \\
& =  \sum_{\gamma \in \langle g 
\rangle} \sum_{\gamma^{\prime \prime} \in \Gamma} \int_{\cal F} \int_{\cal F}  
b(m^\prime (\gamma^{\prime \prime})^{-1} \gamma \gamma^{\prime \prime}, 
m \gamma^{\prime \prime}) \:
a(m \gamma^{\prime \prime}, m^\prime)
\: d\vol(m^\prime)\: d\vol(m) \notag \\
& =  \sum_{\gamma \in \langle g 
\rangle} \sum_{\gamma^{\prime \prime} \in \Gamma} \int_{\cal F} \int_{\cal F}  
b(m^\prime \gamma, m \gamma^{\prime \prime}) \: 
a(m \gamma^{\prime \prime}, m^\prime)
\: d\vol(m) \: d\vol(m^\prime) \notag \\
& =  \sum_{\gamma \in \langle g 
\rangle}  \int_{\cal F} \int_{\hM}  
b(m^\prime \gamma, m) \:
a(m, m^\prime)
\: d\vol(m) \: d\vol(m^\prime) \notag \\
& = \Tr_{\langle g \rangle} (ba). \notag
\end{align} 
It is not hard to justify the steps in (\ref{formally}).
\end{pf}

We will need two slight extensions of the algebra ${\cal A}$. First,
let $E$ be a Hermitian vector bundle on $M$. Put $\hE = \pi^*E$.
For $i \in \{1,2\}$, let $\pi_i$ be projection from $\hM \times \hM$ onto the
$i$-th factor of $\hM$. Let ${\cal A}$ be the convolution algebra of
elements $a \in C^\infty(\hM \times \hM ; \pi_1^* \hE \otimes \pi_2^* \hE^*)$
satisfying the two conditions of Definition \ref{conv}. Equation (\ref{A})
now has to be interpreted as
\begin{equation}
A_{\langle g \rangle}(m) = \sum_{\gamma \in \langle g \rangle} 
\tr \left( (R_\gamma^* a)(m, m) \right) = \sum_{\gamma \in \langle g \rangle} 
\tr \left(a(m \gamma, m) \right).
\end{equation}
Then the proof of Lemma \ref{comm} extends.
Next, we can replace condition 2. of Definition \ref{conv} by the
weaker assumption that
\begin{equation} \label{l1}
\text{for all $R > 0$, } \sup_{d(m, m^\prime) \le R} \sum_{\gamma}
|a(m\gamma, m^\prime)| < \infty.
\end{equation}
Equation (\ref{l1}) is essentially an $l^1$-condition on $a$ with
respect to $\Gamma$. Then we again have a convolution algebra and 
the proof of Lemma \ref{comm} still goes through.

Let $\htriangle_p$ be the $p$-form Laplacian
on $\widehat{M}$. For $t > 0$, let $e^{-t \htriangle_p}$ be the corresponding
heat operator. It has a Schwartz kernel $e^{-t \htriangle_p}(m,m^\prime) \in
\Lambda^p(T_m^* \widehat{M}) \otimes 
(\Lambda^p(T_{m^\prime}^* \widehat{M}))^*$.
By finite propagation speed estimates \cite{Cheeger-Gromov-Taylor (1982)},
$e^{-t \htriangle_p}(m,m^\prime)$ satisfies (\ref{l1}). 

\begin{definition}
Take $g \ne e$. When the limit exists, we define the 
$p$-th delocalized $L^2$-Betti number of $M$ by
\begin{equation}
b_{p, \langle g \rangle}(M) = \lim_{t \rightarrow \infty}
\Tr_{\langle g \rangle} \left( e^{-t \htriangle_p} \right)
\end{equation}
\end{definition}

If we were to take $g = e$ then $b_{p, \langle e \rangle}(M)$ would be the
same as the $p$-th $L^2$-Betti number of $M$ \cite{Atiyah (1976)}.

Let $\{ds^2(u)\}_{u \in [-1,1]}$ 
be a smooth $1$-parameter family of Riemannian metrics
on $M$.  Let $\{*(u)\}_{u \in [-1,1]}$ be the corresponding 
$1$-parameter family of Hodge duality operators.  Then
\begin{align}
\frac{d}{du} \Tr_{\langle g \rangle} \left( e^{-t \htriangle_p} \right) & =
- t \: \Tr_{\langle g \rangle} \left( e^{-t \htriangle_p} \: 
\frac{d\htriangle_p}{du}\right) \\
& = - t \: \Tr_{\langle g \rangle} \left( e^{-t \htriangle_p} \: 
\left(d \left[d^*, *^{-1} \frac{d*}{du}\right] +
\left[d^*, *^{-1} \frac{d*}{du}\right] d \right) \right) \notag \\
& = - 2 t \: \Tr_{\langle g \rangle} \left( e^{-t \htriangle_p} \: 
\left( d d^* - d^* d \right) \: *^{-1} \frac{d*}{du} \right). \notag
\end{align}
Let $\Pi_{Ker(\htriangle_p)}$ be projection onto the $L^2$-kernel of
$\htriangle_p$.
As $\lim_{t \rightarrow \infty} e^{-t \htriangle_p} = \Pi_{Ker(\htriangle_p)}$
and
\begin{equation}
\Pi_{Ker(\htriangle_p)} d d^*  = \Pi_{Ker(\htriangle_p)} d^* d = 0,
\end{equation}
we expect that
\begin{equation} \label{bind}
\lim_{t \rightarrow \infty}
- 2 t \: \Tr_{\langle g \rangle} \left( e^{-t \htriangle_p} \: 
\left( d d^* - d^* d \right) \: *^{-1} \frac{d*}{du} \right) = 0.
\end{equation}
In summary, we have shown the following result.
\begin{proposition}
If (\ref{bind}) can be justified, uniformly in $u$, then
$b_{p,\langle g \rangle}(M)$ is metric-independent and hence a
(smooth) topological invariant of $M$.
\end{proposition}

For the moment, let us assume that $b_{p,\langle g \rangle}(M)$ is 
metric-independent. For $t > 0$, put
\begin{equation} \label{T(t)}
{\cal T}_{\langle g \rangle}(t)  = 
\sum_{p=0}^{d} (-1)^p \: p \: \Tr_{\langle g \rangle} 
\left( e^{-t \htriangle_p} \right).
\end{equation}
Put
\begin{equation} \label{T(infinity)}
{\cal T}_{\langle g \rangle}(\infty) = 
 \sum_{p=0}^{d} (-1)^p \: p \: b_{p,\langle g \rangle}(M). 
\end{equation}
\begin{definition}
Take $g \ne e$.
When the integral makes sense, we define the delocalized $L^2$-analytic
torsion by
\begin{equation} \label{torint}
{\cal T}_{\langle g \rangle}(M) = 
- \int_0^\infty \left( {\cal T}_{\langle g \rangle}(t) - (1-e^{-t})
{\cal T}_{\langle g \rangle}(\infty) \right) \frac{dt}{t}.
\end{equation}
\end{definition}

If we were to take $g = e$ then $T_{\langle g \rangle}(M)$ would formally
be the same, up to a sign, as
the $L^2$-analytic torsion of \cite{Lott (1992b),Mathai (1992)}. 
(The latter requires a zeta-function 
regularization in its definition, but this is not
necessary for the delocalized $L^2$-analytic torsion.) 
It follows from finite propagation
speed arguments that the integrand in (\ref{torint}) is integrable for
small $t$. Thus the question of whether the integral makes sense refers
to large-$t$ integrability.

Let $\{ds^2(u)\}_{u \in [-1,1]}$ 
be a smooth $1$-parameter family of Riemannian metrics
on $M$. Then for any $t > 0$,
as in the proof of \cite[Lemma 8]{Lott (1992b)}, we have
\begin{equation} \label{eqstart}
\frac{d}{du} {\cal T}_{\langle g \rangle}(t)  =
t \frac{d}{dt}
\sum_{p=0}^{d} (-1)^p \:
\Tr_{\langle g \rangle} \left( e^{-t \htriangle_p} *^{-1} \frac{d*}{du} 
\right).
\end{equation}  
Proceeding formally,
\begin{equation}
\frac{d}{du} {\cal T}_{\langle g \rangle}(M) =
\left( \lim_{t \rightarrow 0} - \lim_{t \rightarrow \infty} \right)
\sum_{p=0}^{d} (-1)^p \: 
\Tr_{\langle g \rangle} \left( e^{-t \htriangle_p} *^{-1} \frac{d*}{du} 
\right).
\end{equation}  
As we are assuming that $g \ne e$, it follows again from finite propagation
speed estimates that
\begin{equation}
\lim_{t \rightarrow 0}
\sum_{p=0}^{d} (-1)^p \:  
\Tr_{\langle g \rangle} \left( e^{-t \htriangle_p} *^{-1} \frac{d*}{du} 
\right) = 0.
\end{equation}
When it can justified, we expect that
\begin{equation} \label{largetvar}
\lim_{t \rightarrow \infty}
\sum_{p=0}^{d} (-1)^p \:  
\Tr_{\langle g \rangle} \left( e^{-t \htriangle_p} *^{-1} \frac{d*}{du} 
\right) =
\sum_{p=0}^{d} (-1)^p \:  
\Tr_{\langle g \rangle} \left( \Pi_{Ker(\htriangle_p)} *^{-1} \frac{d*}{du} 
\right).
\end{equation}
Now $\Ker(\htriangle_p)$ can be identified with the $p$-dimensional
(reduced) $L^2$-cohomology group of $M$ and is topological in nature
\cite{Atiyah (1976)}.
In summary, we have shown the following result.
\begin{proposition}
If (\ref{largetvar}) can be justified, uniformly in $u$, 
and if $M$ has vanishing $L^2$-cohomology
then for all $g \ne e$, ${\cal T}_{\langle g \rangle}(M)$ is a (smooth)
topological invariant of $M$.
\end{proposition}

Now let $\hE$ be a $\Gamma$-invariant Clifford module on $\hM$. For
simplicity, we assume that $M$ is spin, with $S$ denoting the spinor bundle,
and that there is a Hermitian vector bundle $V$ on $M$ so that 
$\hE = \pi^*(S \otimes V)$; the general case is similar.
Let $\nabla^S$ be the connection on $S$ coming from the
Levi-Civita connection and let $\nabla^V$ be a Hermitian connection on $V$. 
Let $D$ be the corresponding self-adjoint
Dirac-type operator on sections of $S \otimes V$ and let 
$\hD$ be the lifted operator on sections of $\hE$.  
For $t > 0$, let $e^{-t \hD^2}$ be the corresponding
heat operator.  It has a Schwartz kernel $e^{-t \hD^2}(m,m^\prime) \in
\hE_m \otimes \hE_{m^\prime}^*$. 
Again, using finite propagation estimates it is not hard to see that 
$e^{-t \hD^2}$ satisfies (\ref{l1}).
Given a conjugacy class $\langle g \rangle$ in
$\Gamma$, for $s > 0$ put
\begin{equation} \label{eta(s)}
{\eta}_{\langle g \rangle}(s)  = 
\Tr_{\langle g \rangle} \left( \hD e^{-s^2 \hD^2} \right).
\end{equation}

\begin{definition}
Take $g \ne e$. When the integral makes sense, we define the delocalized
$L^2$-eta invariant by
\begin{equation} \label{etaint}
{\eta}_{\langle g \rangle}(M)  =  \frac{2}{\sqrt{\pi}}
\int_0^\infty {\eta}_{\langle g \rangle}(s) ds.
\end{equation}
\end{definition}
If we were to take $g = e$ then ${\eta}_{\langle g \rangle}(M)$ would be
the same as the $L^2$-eta invariant of Cheeger and Gromov 
\cite{Cheeger-Gromov (1985)}. In this case, it is known that the integral
in (\ref{etaint}) makes sense 
\cite{Bismut-Freed (1986),Peric (1992),Ramachandran (1993)}.
If $g \ne e$ then finite propagation speed arguments show that the
integrand in (\ref{etaint}) is integrable for small-$s$. Thus the question
of whether the integral makes sense refers to large-$s$ integrability.
Equation (\ref{etaint}) was first
considered, when $\Gamma$ is virtually nilpotent,
in \cite[Eqn. (69)]{Lott (1992a)}.

Let $\{ds^2(u)\}_{u \in [-1,1]}$, $\{h^V(u)\}_{u \in [-1,1]}$ and
$\{\nabla^V(u)\}_{u \in [-1,1]}$ 
be smooth $1$-parameter families of
Riemannian metrics on $M$, Hermitian metrics on $V$ and
compatible Hermitian connections on $V$, respectively.
Then for any $s > 0$, one can check that
\begin{equation} \label{eqstart2}
\frac{d}{du} {\eta}_{\langle g \rangle}(s)  =
\frac{d}{ds}
\Tr_{\langle g \rangle} \left( s \: \frac{d\hD}{du} e^{-s^2 \hD^2}
\right).
\end{equation}  
Proceeding formally,
\begin{equation}
\frac{d}{du} {\eta}_{\langle g \rangle}(M) =
\frac{2}{\sqrt{\pi}} \left( \lim_{s \rightarrow \infty} - 
\lim_{s \rightarrow 0} \right)
\Tr_{\langle g \rangle} \left( s \: \frac{d\hD}{du} e^{-s^2 \hD^2}
\right).
\end{equation}  
As we are assuming that $g \ne e$, it follows from finite propagation
speed estimates that
\begin{equation} \label{eqend2}
\frac{2}{\sqrt{\pi}} \lim_{s \rightarrow 0}
\Tr_{\langle g \rangle} \left( s \: \frac{d\hD}{du} e^{-s^2 \hD^2}
\right) = 0.
\end{equation}
In summary, we have shown the following result.
\begin{proposition} \label{etalim}
If 
\begin{equation}
\frac{2}{\sqrt{\pi}} \lim_{s \rightarrow \infty}
\Tr_{\langle g \rangle} \left( s \: \frac{d\hD}{du} e^{-s^2 \hD^2}
\right) = 0
\end{equation}
uniformly in $u$
then $\eta_{\langle g \rangle}(M)$ is independent of $u$.
\end{proposition}

Proceeding very formally,
\begin{equation}
\frac{2}{\sqrt{\pi}} \lim_{s \rightarrow \infty}
\Tr_{\langle g \rangle} \left( s \: \frac{d\hD}{du} e^{-s^2 \hD^2}
\right) = 
2 \Tr_{\langle g \rangle} \left(  \frac{d\hD}{du} \delta(\hD)
\right)
 = 2 \frac{d}{du}
\Tr_{\langle g \rangle} \left( \sign(\hD)
\right).
\end{equation}
Thus we expect that if $\sign(\hD)$ is independent of $u$,
as a $\Gamma$-operator, then $\eta_{\langle g \rangle}(M)$ 
is independent of $u$.
In particular, we expect that this will be true in the following cases :\\
1. If $D$ is the Dirac operator and
for all $u \in [-1,1]$, $(M, ds^2(u))$ has positive scalar curvature.\\
2. If $D$ is the (tangential) signature operator of
\cite{Atiyah-Patodi-Singer (1976)}.

We now give some elementary properties of the delocalized $L^2$-invariants.
\begin{proposition} \label{finite}
Suppose that $\Gamma$ is finite. Let 
$\{\langle g_i \rangle\}$ parametrize the 
conjugacy classes of $\Gamma$.  Let $\rho : \Gamma \rightarrow U(N)$ be
a unitary representation of $\Gamma$. Let $E_\rho$ be the associated flat
Hermitian vector bundle on $M$.
Let $\chi_\rho : \Gamma \rightarrow \C$ be
the character of $\rho$.\\
1. Then
\begin{equation}
b_p(M; E_\rho) = \sum_i \chi_\rho(g_i) \: b_{p, \langle g_i \rangle}(M).
\end{equation}
\noindent
2. Let ${\cal T}(M; E_\rho) \in \R$ be the Ray-Singer 
analytic torsion
\cite{Ray-Singer (1971)}.
Then
\begin{equation}
{\cal T}(M; E_\rho) =
\sum_i \chi_\rho(g_i) \: {\cal T}_{\langle g_i \rangle}(M). 
\end{equation}
\noindent
3. Let $D$ be a Dirac-type operator on $M$. 
Let $\eta(M; E_\rho) \in \R$ be the Atiyah-Patodi-Singer
eta-invariant
\cite{Atiyah-Patodi-Singer (1976)}. Then
\begin{equation}
\eta(M; E_\rho) = \sum_i \chi_\rho(g_i) \: \eta_{\langle g_i \rangle}(M). 
\end{equation} 
\end{proposition}
\begin{pf}
This follows from Fourier analysis on $\Gamma$, as in
\cite[Section 2]{Lott-Rothenberg (1991)}. We omit the details.
\end{pf}

Proposition \ref{finite}.3 shows that when $\Gamma$ is a finite 
group, the delocalized $L^2$-eta invariant has the same information as the
$\rho$-invariant of \cite{Atiyah-Patodi-Singer (1976)}.

\begin{proposition}
1. We have $b_{p, \langle g^{-1} \rangle}(M) = 
\overline{b_{p, \langle g \rangle}(M)}$, 
${\cal T}_{\langle g^{-1} \rangle}(M) = 
\overline{{\cal T}_{\langle g \rangle}(M)}$ and
$\eta_{\langle g^{-1} \rangle}(M) = 
\overline{\eta_{\langle g \rangle}(M)}$. \\
2. Given pairs $(M_1, \Gamma_1)$ and $(M_2, \Gamma_2)$, we have
\begin{equation}
{\cal T}_{\langle g_1, g_2 \rangle}(M_1 \times M_2) = 
\delta_{g_1, e_1} \:  \chi(M_1) \: {\cal T}_{\langle g_2 \rangle}(M_2) +
\delta_{g_2, e_2} \: \chi(M_2) \: {\cal T}_{\langle g_1 \rangle}(M_1)
\end{equation}
and
\begin{align}
{\eta}_{\langle g_1, g_2 \rangle}(M_1 \times M_2) = 
& \delta_{g_1, e_1} \left(  \int_{M_1} \widehat{A}(TM_1) \cup \ch(E_1) \right) 
{\eta}_{\langle g_2 \rangle}(M_2) + \\
& \delta_{g_2, e_2} \left(  \int_{M_2} \widehat{A}(TM_2) \cup \ch(E_2) \right) 
{\eta}_{\langle g_1 \rangle}(M_1). \notag
\end{align}
3. If $d$ is even then ${\cal T}_{\langle g \rangle}(M) = 0$.\\
4. Suppose that $d$ is odd and $D$ is the 
(tangential) signature operator \cite{Atiyah-Patodi-Singer (1976)}.
Then $\eta_{\langle g \rangle}(M) = 0$ if $d \equiv 1 \mod 4$.
\end{proposition}
\begin{pf}
As $e^{-t\triangle_p}$ is self-adjoint and $\Gamma$-invariant,
\begin{align}
\tr \left( e^{-t \htriangle_p} (m\gamma^{-1}, m) \right) = &
\tr \left( e^{-t \htriangle_p} (m, m\gamma) \right) =
\tr \left( e^{-t \htriangle_p} (m\gamma, m) \right)^* \\
= & \overline{\tr \left( e^{-t \htriangle_p} (m\gamma, m) \right)}. \notag
\end{align}
It follows that $b_{p, \langle g^{-1} \rangle}(M) = 
\overline{b_{p,\langle g \rangle}(M)}$. The proof of the rest of 1. is similar.
The proofs of 2., 3. and 4. follow from arguments as in
\cite{Atiyah-Patodi-Singer (1976)} and \cite{Lott (1992b)}.
We omit the details.
\end{pf}
The main results of this paper are given by the following propositions.
\begin{proposition} \label{toshow1}
Suppose that $\langle g \rangle$ is a nontrivial finite conjugacy class. 
Then $b_{p, \langle g \rangle}(M)$ is well-defined and metric-independent. 
If $\Gamma$ is a free abelian group then $b_{p, \langle g \rangle}(M) = 0$.
\end{proposition}

\begin{proposition} \label{toshow2}
Suppose that $\langle g \rangle$ is a nontrivial
finite conjugacy class. Suppose
that $M$ has positive Novikov-Shubin invariants 
\cite{Lott-Lueck (1995),Novikov-Shubin (1986)}.
Then the integrand in (\ref{torint}) is integrable.
If $M$ has vanishing $L^2$-cohomology groups then 
${\cal T}_{\langle g \rangle}(M)$ is metric-independent.
\end{proposition}
\begin{proposition} \label{toshow3}
Suppose that $\langle g \rangle$ is a nontrivial
finite conjugacy class and $D$ is a
Dirac-type operator. Then the integrand in (\ref{etaint}) is integrable.
If $D$ is the (tangential)
signature operator then $\eta_{\langle g \rangle}(M)$ is
metric-independent.  If $D$ is the Dirac operator then 
$\eta_{\langle g \rangle}(M)$ is a locally-constant function on the
space of positive-scalar-curvature metrics on $M$.
\end{proposition}

\begin{proposition} \label{toshow4}
Suppose that $\Gamma$ is virtually nilpotent or 
Gromov-hyperbolic. Suppose that $0 \notin \spec(\htriangle_p)$ or that
$0$ is isolated in $\spec(\htriangle_p)$.
Then for all nontrivial conjugacy classes
$\langle g \rangle$ of $\Gamma$, 
$b_{p, \langle g \rangle}(M)$ is well-defined and metric-independent.
\end{proposition}

\begin{proposition} \label{toshow5}
Suppose that $\Gamma$ is virtually nilpotent or 
Gromov-hyperbolic. Suppose that $0 \notin \spec(\hD)$ or that
$0$ is isolated in $\spec(\hD)$.
Then for all nontrivial conjugacy classes
$\langle g \rangle$ of $\Gamma$, the integrand in (\ref{etaint}) is integrable.
Furthermore, if $D$ is the Dirac operator then
$\eta_{\langle g \rangle}(M)$ is a locally-constant function on the
space of positive-scalar-curvature metrics on $M$.
\end{proposition}

\begin{proposition} \label{toshow6}
Let $M^d$ be a closed oriented hyperbolic manifold.
Let  $\langle g \rangle$ be a nontrivial conjugacy class in $\pi_1(M)$.
Then $b_{p, \langle g \rangle}(M) = 0$ for all $p$.

Suppose that $d = 2n+1$.
Then ${\cal T}_{\langle g \rangle}(M)$ and
$\eta_{\langle g \rangle}(M)$ are well-defined, the latter being
constructed with the (tangential) signature operator.
Let $c$ be the unique closed geodesic in the free homotopy class
specified by $\langle g \rangle$. Let $k \in \Z^+$ be the
multiplicity of $c$, meaning the number of times that $c$ covers a prime
closed geodesic. Let $l$ be the hyperbolic length of $c$. Let
$m \in SO(2n)$ be the parallel transport of a normal vector around $c$. Let
$\sigma_j(m) \in SO(\Lambda^j(\R^{2n}))$ be the action of $m$ on
the exterior power $\Lambda^j(\R^{2n})$.
Then
\begin{equation} \label{hyptor}
{\cal T}_{\langle g \rangle}(M) = \frac{e^{-nl}}{k \det(I-e^{-l} m)} 
\sum_{j=0}^{2n} (-1)^j e^{-l|n-j|} \Tr(\sigma_j(m)). 
\end{equation}
In particular, the marked length spectrum of $M$ can be recovered from
$\{{\cal T}_{\langle g \rangle}(M)\}_{\langle g \rangle \in {\cal C}}$.

If $n$ is even then $\eta_{\langle g \rangle}(M) = 0$. If $n$ is odd,
define angles $\{\theta_j\}_{j=1}^{n}$ by saying that the diagonalization
of $m$ consists of the $2 \times 2$ blocks
\begin{equation}
\left(
\begin{matrix}
\cos \theta_j  & - \sin \theta_j \\
\sin \theta_j & \cos \theta_j
\end{matrix}
\right).
\end{equation}
Put $\mu_j = e^{\frac{l+i\theta_j}{2}}$. 
Then
\begin{equation} \label{hypeta}
\eta_{\langle g \rangle}(M) = \frac{(2i)^{n+1}}{2 \pi k} 
\frac
{\sin(\theta_1) \ldots \sin(\theta_n)}
{|\mu_1 - \mu_1^{-1}|^2 \ldots |\mu_n - \mu_n^{-1}|^2}.
\end{equation}

If $n=1$ then
\begin{equation} \label{n=1}
{\cal T}_{\langle g \rangle}(M) - i \pi \eta_{\langle g \rangle}(M)
= \frac{2}{k} \: \frac{1}{1-\mu_1^2}.
\end{equation}
\end{proposition}

\begin{proposition} \label{toshow7}
Let $Z^n$ be a 
smooth closed even-dimensional manifold and let $\phi$ be a diffeomorphism
of $Z$.  Let
$\phi^*_p \in \Aut(\HH^p(Z; \C))$ be the induced map on cohomology.
Let $M$ be the mapping torus
\begin{equation}
M = (Z \times [0,1])/\{(z,0) \sim (\phi(z),1)\}.
\end{equation}
Put $\Gamma = \Z$, acting on the cyclic cover $\hM$ of $M$.
Define $f : \C \rightarrow \C$ by
\begin{equation}
f(\lambda) = 
\begin{cases}
\lambda &\text{if $|\lambda| \le 1$,} \\
\overline{\lambda^{-1}} &\text{if $|\lambda| > 1$.}
\end{cases}
\end{equation}
Then
\begin{equation} \label{torcoho}
{\cal T}_{\langle k \rangle}(M) = 
\begin{cases}
\frac{1}{k}
\sum_{p=0}^{d} (-1)^p \: \Tr \left[ f(\phi^*_p) \right]^k
&\text{if $k > 0$,}\\
& \\
- \: \frac{1}{k}
\sum_{p=0}^{d} (-1)^p \: \Tr \left[ f(\overline{\phi^*_p}) \right]^{-k}
&\text{if $k < 0$.}
\end{cases}
\end{equation}

Equivalently, let $L(\phi^k)$ be the Lefschetz number of $\phi^k$.
For $z \in \C$ with $|z|$ small, put
$\zeta(z) = \exp \left(\sum_{k > 0} \frac{z^k}{k} \: L(\phi^k) \right)$. 
Then $\zeta(z)$ has a meromorphic continuation
to $z \in \C$, and
\begin{equation} \label{Lefschetz} 
{\cal T}_{\langle k \rangle}(M) =  \int_{S^1} e^{-ik\theta} \: 
\ln |\zeta(e^{i\theta})|^2 \: \frac{d\theta}{2\pi}.
\end{equation}

Suppose that $\phi$ preserves a Dirac-type operator $D_Z$ on $Z$. Let
$D$ be the suspended Dirac-type operator on $M$. Then 
$\eta_{\langle k \rangle}(M)$ is
given in terms of the Atiyah-Bott indices by
\begin{equation}
\eta_{\langle k \rangle}(M) = \frac{i}{k \pi} \: \Tr_s \left( \phi^k 
\Big|_{Ker(D_Z)} \right).
\end{equation}
\end{proposition}

\begin{proposition} \label{toshow8}
Let $Z^n$ be a smooth closed even-dimensional
manifold.  Let $F$ be a finite group and let
$\hZ$ be a connected normal $F$-cover of $Z$. Let $\phi$ be a diffeomorphism of
$Z$ and let $M$ be the mapping torus of $\phi$. Let $\hphi$ be a lift of
$\phi$ to $\hZ$ and let $\alpha \in
\Aut(F)$ be the automorphism defined by
\begin{equation}
\hphi(zf) = \hphi(z) \alpha^{-1}(f)
\end{equation}
for all $z \in \hZ$ and $f \in F$.
Put $\Gamma = F \widetilde{\times}_{\alpha} \Z$, 
acting on $\hZ \times \R$ on the right by
\begin{equation} \label{action}
(z,t) \cdot (f,k) = (\hphi^{k}(zf), t+k).
\end{equation}
For $k \in \Z$, define an equivalence relation $\sim_k$ on $F$ by saying that
$f \sim_{k} f^\prime$ if there exists a $\gamma \in F$ such that
$\gamma f \alpha^{k}(\gamma^{-1}) = f^\prime$. 
Let $[f]_k$ be the equivalence class of $f \in F$ and let
$|[f]_k|$ be its cardinality.
Let $I_k(f) \in \Z$ be the
Nielsen fixed point index of $\phi^k$, evaluated at $[f]_k$.
If $\rho$ is a
finite-dimensional irreducible unitary representation of 
$F \widetilde{\times}_{\alpha} \Z$, let
$\chi_\rho$ be its character. For $z \in \C$ with $|z|$ small, put
\begin{equation}
\zeta_{\rho}(z) = \exp \left(
\sum
\begin{Sb}
f \in F \\
k > 0
\end{Sb}
\frac{z^k}{k} \: \chi_\rho (f,k) \: \frac{I_{k}(f)}{|[f]_k|} \right).
\end{equation}
Then $\zeta_{\rho}(z)$ has a meromorphic continuation to $z \in \C$ 
and
\begin{equation} \label{eq8}
\sum_{f,k} \chi_\rho(f,k) \: {\cal T}_{\langle f,k \rangle}(M) =
\ln |\zeta_{\rho}(1)|^2.
\end{equation}
Knowing (\ref{eq8}) for all $\rho$ 
determines $\{ {\cal T}_{\langle f,k \rangle}(M) 
\}_{(f,k) \in \Gamma}$.
\end{proposition}

Before proceeding with the proofs, 
let us mention why the existence problem for the delocalized $L^2$-invariants
is more difficult than for the ordinary $L^2$-invariants.
The algebraic origin of the problem is as follows.
Let $\Gamma$ be a countable discrete group and
consider the group algebra $\C \Gamma$ of finite sums $\sum_{\gamma \in \Gamma}
c_\gamma \gamma$, with $c_\gamma \in \C$. Define an involution on
$\C \Gamma$ by 
\begin{equation}
\left(\sum_{\gamma \in \Gamma} c_\gamma \gamma \right)^* =
\sum_{\gamma \in \Gamma} \overline{c_\gamma} \: \gamma^{-1}. 
\end{equation}
If $\langle g \rangle$ is a conjugacy class in $\Gamma$, the
linear functional $\tau_{\langle g \rangle} : \C \Gamma \rightarrow \C$ given
by 
\begin{equation} \label{tauu}
\tau_{\langle g \rangle}\left( \sum_{\gamma \in \Gamma} c_\gamma \gamma \right)
= \sum_{\gamma \in \langle g \rangle} c_\gamma
\end{equation}
satisfies $\tau_{\langle g \rangle}(ab) = \tau_{\langle g \rangle}(ba)$ 
for all $a,b \in \C \Gamma$. Furthermore, if $g = e$ then
$\tau_{\langle e \rangle}(a^* a) \ge 0$ and $\tau_{\langle e \rangle}$
extends to a continuous linear functional on the group von Neumann algebra.
These last two properties of $\tau_{\langle e \rangle}$, 
which are crucial for the usual $L^2$-invariants, generally fail for
$\tau_{\langle g \rangle}$ if $g \ne e$.

\section{Proofs of Propositions \ref{toshow1}, \ref{toshow2} and \ref{toshow3}}
\label{proof1}
Let $\langle g \rangle$ be a finite conjugacy class in $\Gamma$.
Put
\begin{equation}
A = \sum_{\gamma \in \langle g \rangle} R_{\gamma}^*.
\end{equation}
Then $A$ is a bounded operator on $\Omega^*(\hM)$ which commutes with
$R_{\gamma^\prime}^*$ for all $\gamma^\prime \in \Gamma$. Letting
$\Tr_\Gamma$ denote the $II_\infty$-trace \cite{Atiyah (1976)}, we have
\begin{equation}
\Tr_{\langle g \rangle} \left( e^{-t \htriangle_p} \right) =
\Tr_\Gamma \left( A e^{-t \htriangle_p} \right).
\end{equation}
Thus
\begin{equation}
b_{p,\langle g \rangle}(M) = 
\Tr_\Gamma \left( A \Pi_{Ker(\htriangle_p)} \right). 
\end{equation}
Hence $b_{p,\langle g \rangle}(M)$ is well-defined. It follows by standard
arguments that it is metric-independent.

Similarly,
\begin{equation}
{\cal T}_{\langle g \rangle}(t) = \sum_{p=0}^n (-1)^p \: p \:
\Tr_\Gamma \left( A e^{-t \htriangle_p} \right). 
\end{equation}
Then the integrand of (\ref{torint}) is
\begin{equation}
\frac{1}{t} \sum_{p=0}^n (-1)^p \: p \:
\Tr_\Gamma \left[ A \left( e^{-t \htriangle_p} - (1-e^{-t})
\Pi_{Ker(\htriangle_p)} \right) \right].
\end{equation}
Now 
\begin{equation}
\Big| \Tr_\Gamma \left[ A \left( e^{-t \htriangle_p} - 
\Pi_{Ker(\htriangle_p)} \right) \right] \Big| \: \le \:
\parallel A \parallel \Tr_\Gamma \left[ \left( e^{-t \htriangle_p} - 
\Pi_{Ker(\htriangle_p)} \right) \right].
\end{equation}
By assumption, there is an $\alpha_p > 0$ such that for large $t$,
\begin{equation}
\Tr_\Gamma \left[ e^{-t \htriangle_p} - 
\Pi_{Ker(\htriangle_p)} \right] \: \le \: t^{-\alpha_p/2}.
\end{equation}
It follows that the integrand in (\ref{torint}) is integrable. If $M$ has
vanishing $L^2$-cohomology groups then it follows as in \cite{Lott (1992b)}
that ${\cal T}_{\langle g \rangle}(M)$ is metric-independent.

Now let $D$ be a Dirac-type operator.  Then we have

\begin{equation} 
{\eta}_{\langle g \rangle}(s)  = 
\Tr_\Gamma \left( A \hD e^{-s^2 \hD^2} \right)
\end{equation}
Thus
\begin{equation}
\Big| {\eta}_{\langle g \rangle}(s) \Big| \: \le \:
\parallel A \parallel \Tr_\Gamma \left(  \big| \hD \big| 
e^{-s^2 \hD^2} \right).
\end{equation}
It follows as in \cite[Theorem 3.1.1]{Ramachandran (1993)} 
that the integrand in
(\ref{etaint}) is integrable. The rest of Proposition \ref{toshow3} follows
as in \cite{Cheeger-Gromov (1985)}. 

Finally, to finish the proof of Proposition \ref{toshow1},
suppose that $\Gamma = \Z^l$. Let $\hGamma$ be the Pontryagin dual of
$\Gamma$. Then an element $\rho_\theta$ of $\hGamma$ is a representation
$\rho_\theta : \Gamma \rightarrow U(1)$ of the form
\begin{equation}
\rho_\theta (\Vec{k}) = e^{i \Vec{k} \cdot \Vec{\theta}}.
\end{equation}
Let $E_\theta$ be the associated flat unitary line bundle on $M$ and
let $\triangle_{p,\theta}$ be the Laplacian on $\Omega^p(M; E_\theta)$. 
It follows as in \cite[Proposition 38]{Lott (1992b)} that
\begin{equation}
\Tr_{\langle \Vec{m} \rangle} \left( e^{-t \htriangle_p} \right) = 
\int_{\hGamma} e^{-i \Vec{m} \cdot \Vec{\theta}} \:
\Tr \left( e^{-t \triangle_{p,\theta}} \right) \: 
\frac{d^l\theta}{(2\pi)^l}.
\end{equation}
From \cite[Chapter XII]{Reed-Simon (1980)}, the eigenvalues 
$\lambda_i(\theta)$ form a sequence of nonnegative algebraic
functions locally on $\hGamma$. (This corrects a claim in 
\cite{Lott (1992b)} that they are local analytic functions, which
is only guaranteed if $l = 1$.) It follows that there is
convergence in $L^1(\hGamma)$ :
\begin{equation}
\lim_{t \rightarrow \infty} \Tr \left( e^{-t \triangle_{p,\theta}} \right)
= b_p^{(2)}(M),
\end{equation}
where $b_p^{(2)}(M)$ is the number of such algebraic
functions which equal the zero function.
Hence if $\Vec{m} \ne 0$ then $b_{p, \langle \Vec{m} \rangle}(M) = 0$.
\qed

\section{Proofs of Propositions \ref{toshow4} and \ref{toshow5}}

Let $\Lambda$ be the reduced group $C^*$-algebra of $\Gamma$.
We assume that there is an algebra ${\frak B}$ such that\\
1. $\C \Gamma \subseteq {\frak B} \subseteq \Lambda$.\\
2. ${\frak B}$ is the projective limit of a sequence
\begin{equation}
\ldots \rightarrow B_{j+1} \rightarrow B_j \rightarrow \ldots \rightarrow
B_0
\end{equation}
of Banach algebras $(B_j, |\cdot|_j)$ with unit, and $B_0 = \Lambda$.\\
3. If $i_j : {\frak B} \rightarrow B_j$ is the induced homomorphism
then $i_0$ is injective with dense image and ${\frak B}$ is closed under
the holomorphic functional calculus in $\Lambda$.\\
4. Given a conjugacy class $\langle g \rangle$ of $\Gamma$, define
$\tau_{\langle g \rangle} : \C \Gamma \rightarrow \C$
as in (\ref{tauu}).
Then for all $\langle g \rangle \in {\cal C}$, $\tau_{\langle g \rangle}$
extends to a continuous linear functional on ${\frak B}$, which we
again denote by $\tau_{\langle g \rangle}$.

The topology on ${\frak B}$ comes from the submultiplicative
seminorms $\parallel \cdot \parallel_j = |i_j(\cdot)|_j$.
Condition 4 is equivalent to saying that $\HC^0(\C \Gamma) = 
\HC^0({\frak B})$. If $\Gamma$ is a virtually nilpotent or 
Gromov-hyperbolic group then conditions 1-4 are known to be satisfied
by the rapid-decay algebra ${\frak B}$ \cite[p. 397]{Ji (1993)}.

If ${\frak E}$ is a finitely-generated right projective ${\frak B}$-module
then there is a continuous trace
\begin{equation} \label{contrace}
\Tr : \End_{\frak B}({\frak E}) \rightarrow {\frak B}/
\overline{[{\frak B},{\frak B}]}.
\end{equation}
Explicitly, suppose that ${\frak E} = {\frak B}^N e$ for some $N > 0$ and
some projection
$e \in M_N({\frak B})$. If $A \in \End_{\frak B}({\frak E})$, we can
think of $A$ as an element of $M_N({\frak B})$ satisfying
$A = eA = Ae$. Then
\begin{equation}
\Tr(A) = \sum_{i=1}^N A_{ii} \mod \overline{[{\frak B},{\frak B}]}.
\end{equation}
(We quotient by the closure of $[{\frak B},{\frak B}]$ to ensure that the
trace takes value in a Fr\'echet space.) 

As $\Lambda$ is a $C^*$-algebra, there is a calculus of
$\Lambda$-pseudodifferential operators  on $M$
\cite{Mischenko-Fomenko (1979)}. 
Suppose that $E^1$ is a smooth
$\Lambda$-vector bundle on $M$, meaning the fibers of $E^1$ are all
isomorphic to a fixed finitely-generated projective right $\Lambda$-module
${\frak E}^1$ and the transition functions are smooth functions with
value in $\Aut_{\Lambda}
({\frak E}^1)$. Let $E^2$ be another smooth $\Lambda$-vector bundle on $M$.
The elements of the pseudodifferential algebra 
$\Psi^\infty_{\Lambda}(M; E^1, E^2)$ map smooth 
sections of $E^1$ to smooth sections of
$E^2$ and commute with the $\Lambda$-action.

In \cite[Section 6.1]{Lott (1996)} we extended this to
a calculus of ${\frak B}$-pseudodifferential operators and proved some
basic properties of such operators.
We only state the necessary facts, referring to  \cite{Lott (1996)} for
details.

Let ${\cal E}^1$ and ${\cal E}^2$ be smooth
${\frak B}$-vector bundles on $M$. By an extension of the Serre-Swan theorem,
we can write ${\cal E}^1 = 
(M \times {\frak B}^N) e^1$ for some $N > 0$ and
some projection $e^1 \in C^\infty(M; M_N({\frak B}))$.
Define a $B_j$-vector bundle by
$E^1_j = (M \times B_j^N) i_j(e^1)$. Then ${\cal E}^1$ is the projective 
limit of
\begin{equation}
\ldots \rightarrow E^1_{j+1} \rightarrow E^1_j \rightarrow \ldots \rightarrow
E^1_0,
\end{equation}
and similarly for ${\cal E}^2$.

For each $j \ge 0$, there is an algebra $\Psi^\infty_{B_j}(M; E^1_j, E^2_j)$ 
of $B_j$-pseudodifferential operators.  The algebra 
$\Psi^\infty_{\frak B}(M; {\cal E}^1, {\cal E}^1)$ 
of ${\frak B}$-pseudodifferential operators is the projective limit of
\begin{equation}
\ldots \rightarrow \Psi^\infty_{B_j}(M; E^1_{j+1}, E^2_{j+1}) \rightarrow 
\Psi^\infty_{B_j}(M; E^1_j, E^2_j) \rightarrow \ldots \rightarrow
\Psi^\infty_{B_0}(M; E^1_0, E^2_0).
\end{equation}

Let ${\cal E}$ be a ${\frak B}$-vector bundle on $M$. Given $T \in
\Psi^\infty_{\frak B}(M; {\cal E}, {\cal E})$, let $i_j(T)$ be its
image in $\Psi^\infty_{B_j}(M; E_j, E_j)$.

\begin{proposition} \label{spec}
\cite[Proposition 19]{Lott (1996)} If $i_0(T)$ is invertible
in $\Psi^\infty_{B_0}(M; E_0, E_0)$ then $T$ is invertible in
$\Psi^\infty_{\frak B}(M; {\cal E}, {\cal E})$.
\end{proposition}

Note that $\Psi^{-\infty}_{\frak B}(M; {\cal E}, {\cal E})$ is an algebra in
its own right (without unit) of smoothing operators. 
Given $T \in \Psi^{-\infty}_{\frak B}(M; 
{\cal E}, {\cal E})$, let $\sigma_{\Psi^{-\infty}}(T)$ denote its spectrum
in $\Psi^{-\infty}_{\frak B}(M; {\cal E}, {\cal E})$ and let
$\sigma_{\Psi^{\infty}}(T)$ denote its spectrum
in $\Psi^{\infty}_{\frak B}(M; {\cal E}, {\cal E})$.

\begin{lemma} \cite[Lemma 2]{Lott (1996)}
$\sigma_{\Psi^{-\infty}}(T) = \sigma_{\Psi^{\infty}}(T)$.
\end{lemma}

Consider the algebra ${\frak A}$ of integral operators whose kernels
$K(m_1, m_2) \in \Hom_{\frak B}({\cal E}_{m_2},{\cal E}_{m_1})$ are 
continuous in $m_1$ and $m_2$, with multiplication
\begin{equation}
(K K^\prime)(m_1, m_2) = \int_Z K(m_1, m) K^\prime(m, m_2) \: d\vol(m).
\end{equation}
Let $A_j$ be the analogous algebra with continuous kernels
$K(m_1, m_2) \in \Hom_{B_j}((E_j)_{m_2},(E_j)_{m_1})$. Give
$\Hom_{B_j}((E_j)_{m_2},(E_j)_{m_1})$ the Banach space norm $| \cdot |_j$ 
induced from $\Hom(B_j^N, B_j^N)$. Define a norm $| \cdot |_j$ on $A_j$ by
\begin{equation}
|K|_j = (\vol(M))^{-1} \max_{m_1, m_2 \in M} |K(m_1, m_2)|_j.
\end{equation}
Then one can check that $A_j$ is a Banach algebra (without unit).
Furthermore, ${\frak A}$ is the projective limit of $\{ A_j \}_{j \ge 0}$.
Any smoothing operator
$T \in \Psi^{-\infty}_{\frak B}(M; {\cal E}, {\cal E})$ gives an
element of ${\frak A}$ through its Schwartz kernel. Let
$\sigma_{\frak A}(T)$ be its spectrum in ${\frak A}$.

\begin{lemma} \cite[Lemma 3]{Lott (1996)}
$\sigma_{\frak A}(T) = \sigma_{\Psi^{-\infty}}(T)$.
\end{lemma}

Define a continuous trace $\TR : {\cal A} \rightarrow {\frak B}/
\overline{[{\frak B},{\frak B}]}$ by
\begin{equation}
\TR(K) = \int_M \Tr(K(m,m)) \: d\vol(m).
\end{equation}

Suppose that $0 \notin \spec(\htriangle_p)$ or that
$0$ is isolated in
$\spec(\htriangle_p)$. Let ${\cal D}$ be the ${\frak B}$-vector bundle 
$\hM \times_\Gamma {\frak B}$ on $M$ and put ${\cal E} = 
\Lambda^p(T^*M) \otimes {\cal D}$.
We can lift $\triangle_p$ from $M$ to a differential operator
$\widetilde{\triangle}_p \in \Psi^2_{\frak B}(M; {\cal E}, {\cal E})$.
Then for all $t > 0$, $e^{- t \widetilde{\triangle}_p} \in 
\Psi^{-\infty}_{\frak B}(M; {\cal E}, {\cal E})$.
As in \cite[Section 3]{Lott (1992)}, one can show that
\begin{equation}
\Tr_{\langle g \rangle} \left( e^{- t \widehat{\triangle}_p} \right)
= \tau_{\langle g \rangle} \left( \TR \left(
e^{- t \widetilde{\triangle}_p} \right) \right).
\end{equation}

Let $E$ be the analogous $\Lambda$-vector bundle on $M$ whose fiber over
$m \in M$ is isomorphic to $\Lambda^p(T_m^*M) \otimes \Lambda$. 
Recall that $i_0(\widetilde{\triangle}_p)$ is the extension of 
$\widetilde{\triangle}_p$ to an element of $\Psi^2_{\Lambda}(M; E, E)$.
Let $N(\Gamma)$ denote the group von Neumann algebra of $\Gamma$ 
\cite{Atiyah (1976)}. 
Let $\overline{E}$ be the natural $N(\Gamma)$-vector bundle on 
$M$ whose fiber over
$m \in M$ is isomorphic to $\Lambda^p(T_m^*M) \otimes N(\Gamma)$.
Let $\overline{\triangle}_p$ be the
extension of $\widetilde{\triangle}_p$ to an element of
$\Psi^2_{N(\Gamma)}(M; \overline{E}, \overline{E})$. As $L^2(\hM)$
is isomorphic to the $L^2$-sections of the Hilbert bundle
$\hM \times_\Gamma l^2(\Gamma)$, it follows that $\sigma(\htriangle_p) =
\sigma(\overline{\triangle}_p)$. 
As the $C^*$-algebra $\Lambda$ is a closed subalgebra
of $N(\Gamma)$, it follows that 
$\sigma(\overline{\triangle}_p) = \sigma(i_0(\widetilde{\triangle}_p))$.
Using Proposition \ref{spec}, it now follows that 
$0 \notin \sigma(\widetilde{\triangle_p})$ or that $0$ is isolated in 
$\sigma(\widetilde{\triangle_p})$. Let $c$ be a small loop around $0 \in \C$,
oriented counterclockwise. The projection onto 
$\Ker(\widetilde{\triangle}_p)$ is
\begin{equation}
\Pi_{Ker(\widetilde{\triangle}_p)} = 
\frac{1}{2\pi i} \int_c \frac{dz}{z-
\widetilde{\triangle}_p}.
\end{equation}
It follows from arguments as in \cite{Mischenko-Fomenko (1979)} that
$\Ker(\widetilde{\triangle}_p)$ is a finitely-generated projective
right ${\frak B}$-module. Let $\widetilde{\triangle}_p^\prime$ be the
compression of $\widetilde{\triangle}_p$ onto 
$\Image \left(I-\Pi_{Ker(\widetilde{\triangle}_p)} \right)$.
Then in terms of the trace of (\ref{contrace}),
\begin{equation}
\Tr_{\langle g \rangle} \left( e^{- t \widehat{\triangle}_p} \right)
= \tau_{\langle g \rangle}
\left( \Tr \left( I_{Ker(\widetilde{\triangle}_p)} \right) \right) +
\tau_{\langle g \rangle} 
\left( \TR \left( e^{-t \widetilde{\triangle}_p^\prime} \right) \right).
\end{equation}
As $\tau$ and $\TR$ are continuous, it suffices to show that
there is some $j$ such that the $A_j$-norm 
of $e^{-t \widetilde{\triangle}_p^\prime}$ is rapidly-decreasing
in $t$. 

As $\{ e^{- t \widetilde{\triangle}_p^\prime} \}_{t > 0}$ gives a 
$1$-parameter semigroup in the Banach algebra $A_j$, it follows from
\cite[Theorem 1.22]{Davies (1980)} that the number
\begin{equation}
a = \lim_{t \rightarrow \infty} t^{-1} \ln \Big|
e^{-t \widetilde{\triangle}_p^\prime}
\Big|_j
\end{equation}
exists.  Furthermore, for all $t > 0$, the spectral radius of
$e^{-t \widetilde{\triangle}_p^\prime}$ is $e^{at}$. Let $\lambda_0 > 0$ be 
the infimum of the spectrum of the generator
$\widetilde{\triangle}_p^\prime$. Then by the spectral mapping
theorem, the spectral radius of
$e^{-t \widetilde{\triangle}_p^\prime}$ is $e^{-t \lambda_0}$. Thus
there is a constant $C > 0$
such that for $t > 1$,
\begin{equation} \label{decay}
\Big| e^{-t \widetilde{\triangle}_p^\prime} 
\Big|_j \le C \: e^{- \lambda_0 t/2}.
\end{equation}
Hence 
\begin{equation}
b_{p, \langle g \rangle}(M) = \tau_{\langle g \rangle}
\left( \Tr \left( I_{Ker(\widetilde{\triangle}_p)} \right) \right)
\end{equation}
is well-defined. By similar arguments one can justify (\ref{bind}), showing
that $b_{p, \langle g \rangle}(M)$ is metric-independent.

Now let $D$ be a Dirac-type operator on $M$ such that $0 \notin \spec(\hD)$ or
$0$ is isolated in
$\spec(\hD)$. Put ${\cal E} = 
S \otimes V \otimes {\cal D}$.
We can lift $D$ to a differential operator
$\widetilde{D} \in \Psi^1_{\frak B}(M; {\cal E}, {\cal E})$.
Then for all $s > 0$, $\widetilde{D} e^{-s^2 \widetilde{D}^2} \in 
\Psi^{-\infty}_{\frak B}(M; {\cal E}, {\cal E})$.
As in \cite[Section 3]{Lott (1992)}, one can show that
\begin{equation}
{\eta}_{\langle g \rangle}(s) = \tau_{\langle g \rangle} \left( \TR \left(
\widetilde{D} e^{-s^2 \widetilde{D}^2} \right) \right).
\end{equation}
From finite-propagation estimates, we know that
${\eta}_{\langle g \rangle}(s)$ is integrable for small-$s$. Hence we
must show that $\tau_{\langle g \rangle} \left( \TR \left(
\widetilde{D} e^{-s^2 \widetilde{D}^2} \right) \right)$ is integrable for
large-$s$. It suffices to show that
there is some $j$ such that the $A_j$-norm 
of $\widetilde{D} e^{-s^2 \widetilde{D}^2}$ is rapidly-decreasing
in $s$. 

Let $E$ be the natural $\Lambda$-vector bundle on $M$ whose fiber over
$m \in M$ is isomorphic to $S_m \otimes V_m \otimes \Lambda$. 
Recall that $i_0(\widetilde{D})$ is the extension of 
$\widetilde{D}$ to an element of $\Psi^1_{\Lambda}(M; E, E)$. As before,
$\sigma(\hD) = \sigma(i_0(\widetilde{D}))$.
Using Proposition \ref{spec}, it now follows that 
$0 \notin \sigma(\widetilde{D})$ or that $0$ is isolated in 
$\sigma(\widetilde{D})$. Let $c$ be a small loop around $0 \in \C$,
oriented counterclockwise. The projection on $\Ker(\widetilde{D})$ is
\begin{equation}
\Pi_{Ker(\widetilde{D})} = \frac{1}{2\pi i} \int_c \frac{dz}{z-
\widetilde{D}}.
\end{equation}

Let $\widetilde{D}^\prime$ be the compression of $\widetilde{D}$ onto
$\Image \left(I-\Pi_{Ker(\widetilde{D})} \right)$. Then
\begin{equation}
\widetilde{D} e^{-s^2 \widetilde{D}^2} =
0_{Ker(\widetilde{D})} \oplus
\widetilde{D}^\prime e^{-s^2 \widetilde{D}^{\prime 2}}.
\end{equation}
Hence we may as well assume that $0 \notin \spec(\widetilde{D})$.
Let $\lambda_0 > 0$ be 
the infimum of the spectrum of the generator
$\widetilde{D}^2$.
As in (\ref{decay}), there is a constant $C > 0$
such that for $t > 1$,
\begin{equation}
\Big| e^{- t \widetilde{D}^2} \Big|_j \le C \: e^{- \lambda_0 t/2}.
\end{equation}
As
\begin{equation}
\Big| \widetilde{D} e^{-s^2 \widetilde{D}^2} \Big|_j \le
\Big| \widetilde{D} e^{- \widetilde{D}^2} \Big|_j \cdot
\Big| e^{- (s^2 - 1) \widetilde{D}^2} \Big|_j,
\end{equation}
it follows that $\eta_{\langle g \rangle}(s)$ is large-$s$ integrable.

Suppose that $\{ds^2(u)\}_{u \in [-1,1]}$ is a smooth $1$-parameter of 
positive-scalar-curvature metrics on $M$. Let $D(u)$ be the Dirac operator
on $M$. Then for all $u \in [-1,1]$, $\hD(u)$ is invertible. Using the
above methods, one sees that
\begin{equation}
\lim_{s \rightarrow \infty} s \: \tau_{\langle g \rangle} \left( \TR \left(
\frac{d\widetilde{D}}{du} \: e^{-s^2 \widetilde{D}^2} \right) \right) = 0.
\end{equation} 
From Proposition \ref{etalim}, 
$\eta_{\langle g \rangle}(M)$ is independent of $u$.
\qed \\ \\
\noindent
{\bf Remark : } If $\Gamma$ is virtually nilpotent, one can also prove
Propositions \ref{toshow4} and \ref{toshow5} using finite propagation
speed estimates on $\hM$, as in \cite{Lott (1992a)}. This does
not work when $\Gamma$ is Gromov-hyperbolic, which is why we use the more
indirect method of proof above.

\section{Proof of Proposition \ref{toshow6}}

As in \cite[Section 2]{Fried (1986)}, the 
Selberg trace formula implies that there are functions
$\{G_j(t)\}_{j=0}^{d-1}$ of the form
\begin{equation} \label{Sel1}
G_j(t) = a_j \: t^{-1/2} \: e^{- \frac{l^2}{4t}} \: e^{-tc_j^2}
\end{equation}
so that for $0 \le j \le d$,
\begin{equation} \label{Sel2}
\Tr_{\langle g \rangle} \left( e^{-t \htriangle_j} \right) = 
G_j(t) + G_{j-1}(t).
\end{equation}
Here $a_j$ and $c_j$ are nonnegative constants whose exact values are
not important for the moment. It is clear from (\ref{Sel1}) and (\ref{Sel2})
that $b_{p,\langle g \rangle}(M)$ vanishes for all $p$.

Now suppose that the dimension of $M$ is $d = 2n+1$.
The Selberg trace formula gives the following result.

\begin{proposition} \cite[Theorem 2]{Fried (1986)}
For $0 \le j \le 2n$, put $c_j = |n-j|$ and
\begin{equation}
G_t (\sigma_j) = \frac{\Tr(\sigma_j(m))}{k \: \det(I-e^{-l}m)} \:
\frac{1}{\sqrt{4\pi t}} \: e^{-\frac{l^2}{4t}} \: e^{-tc_j^2} \:
e^{-nl}.
\end{equation}
Then for $0 \le j \le d$,
\begin{equation} 
\Tr_{\langle g \rangle} \left( e^{-t \htriangle_j} \right) = 
G_t(\sigma_j) + G_t(\sigma_{j-1}).
\end{equation}
\end{proposition}

Hence from (\ref{T(t)}), 
\begin{align}
{\cal T}_{\langle g \rangle}(t) & = \sum_{p=0}^d (-1)^j \: j \: \left[
G_t(\sigma_j) + G_t(\sigma_{j-1}) \right] \\
& = \sum_{j=0}^{2n} (-1)^{j+1} \: G_t(\sigma_j). \notag
\end{align}
For $l > 0$,
\begin{equation}
\int_0^\infty \frac{1}{\sqrt{4 \pi t}} \: e^{-\frac{l^2}{4t}} \:
e^{-tc^2} \: \frac{dt}{t} = \frac{e^{-lc}}{l}. 
\end{equation}
Equation (\ref{hyptor}) now follows from (\ref{torint}).

For $r > 0$, we have
\begin{align}
{\cal T}_{\langle g^r \rangle}(M) & = \frac{e^{-nrl}}{k \det(I-e^{-rl} m^r)} 
\sum_{j=0}^{2n} (-1)^j e^{-rl|n-j|} \Tr(\sigma_j(m^r)) \\
& = \frac{(-1)^n}{k} \: e^{-nrl} \: \Tr(\sigma_n(m^r)) + 
O\left( e^{-2nrl} \right). \notag
\end{align}
Then
\begin{equation} \label{asymp}
l = \frac{1}{n} \sup \{ \alpha \in \R : \big|
{\cal T}_{\langle g^r \rangle}(M) \big| = O(e^{-\alpha r}) \}  
\end{equation}
Hence one recovers the marked length spectrum of $M$ from
$\{{\cal T}_{\langle g \rangle}(M)\}_{\langle g \rangle \in {\cal C}}$.

Let $D$ be the tangential signature operator.  By
\cite[Theorem 2.1]{Millson (1978)},
\begin{equation}
\eta_{\langle g \rangle}(s) = (2i)^{n} \frac{2 \pi i}{k} \:
\frac
{l^2 \: \sin(\theta_1) \ldots \sin(\theta_n)}
{|\mu_1 - \mu_1^{-1}|^2 \ldots |\mu_n - \mu_n^{-1}|^2} \:
\frac{e^{-\frac{l^2}{4s^2}}}{(4\pi)^{3/2} s^3}.  
\end{equation}
Equation (\ref{hypeta}) now follows from (\ref{etaint}). 

Equation (\ref{n=1}) comes from a straightforward computation. \qed \\ \\
\noindent
{\bf Remark : } Using the results of Moscovici and Stanton 
\cite{Moscovici-Stanton (1989),Moscovici-Stanton (1991)}, one can 
extend Proposition \ref{toshow6} to general nonpositively-curved
locally symmetric spaces.

\section{Proof of Proposition \ref{toshow7}}

Let $\pi : M \rightarrow S^1$ be the natural projection map.
For $e^{i \theta} \in U(1)$, let $E_\theta$ be the flat complex
line bundle on $S^1$ with holonomy $e^{i \theta}$. Let
$T(\theta) \in \R$ be the Ray-Singer analytic torsion of $M$, computed
with the flat bundle $\pi^*(E_\theta)$. As in 
\cite[Section VI]{Lott (1992b)}, it follows from Fourier analysis that
\begin{equation} \label{eq1} 
{\cal T}_{\langle k \rangle}(M) = \int_{S^1} e^{-ik\theta} \: T(\theta) \:
\frac{d\theta}{2\pi}.
\end{equation}
From \cite[Section 3]{Milnor (1968)} and the Cheeger-M\"uller theorem
\cite{Cheeger (1979),Mueller (1978)},
\begin{equation} \label{eq2}
T(\theta) = \sum_{p=0}^n (-1)^p \: \ln |\det(I- e^{i\theta} \phi_p^*)|^{-2}.
\end{equation}
Given $\lambda \ne 0$, if $k > 0$ then
\begin{equation} \label{eq3}
\int_{S^1} e^{-ik\theta} \:
\ln |1 - e^{i \theta} \lambda|^{-2} \: \frac{d\theta}{2 \pi} =
\frac{f(\lambda^k)}{k}
\end{equation}
and if $k < 0$ then
\begin{equation} \label{eq4}
\int_{S^1} e^{-ik\theta} \:
 \ln |1 - e^{i \theta} \lambda|^{-2} \: \frac{d\theta}{2 \pi} =
- \: \frac{f(\overline{\lambda}^{-k})}{k}.
\end{equation}
Equation (\ref{torcoho}) follows from combining (\ref{eq1})-(\ref{eq4}). 

By standard arguments,
\begin{equation} \label{Milnor}
\zeta(z) = \prod_{p=0}^n  \det \left( I - z \phi_p^* \right)^{(-1)^{p+1}}.
\end{equation}
Equation (\ref{Lefschetz}) follows from (\ref{eq1}),
(\ref{eq2}) and (\ref{Milnor}).

Now suppose that $\phi$ preserves $D_Z$. It follows that $\phi$ is an 
isometry of $Z$ with respect to the Riemannian metric defining $D_Z$.
In terms of the coordinates
$(u,z)$ on $\hM = \R \times Z$, we can write 
\begin{equation}
\hD = 
\begin{pmatrix}
-i \partial_u & D_{Z,-} \\
D_{Z,+} & i \partial_u
\end{pmatrix}.
\end{equation}
Then
\begin{equation}
{\hD}^2 = 
\begin{pmatrix}
- \partial_u^2 + D_{Z,-} D_{Z,+} & 0 \\
0 & - \partial_u^2 + D_{Z,+} D_{Z,-}
\end{pmatrix}
\end{equation}
and 
\begin{align}
& e^{-s^2 {\hD}^2} ((u, z), (u^\prime,z^\prime
)) = \\
& \hspace{.5in} \begin{pmatrix}
\frac{1}{\sqrt{4\pi s^2}} \: e^{-\frac{(u - u^\prime)^2}{4s^2}} \: 
e^{-s^2 D_{Z,-} D_{Z,+}} (z, z^\prime) & 0 \\
0 & \frac{1}{\sqrt{4\pi s^2}} \: e^{-\frac{(u - u^\prime)^2}{4s^2}} \: 
e^{-s^2 D_{Z,+} D_{Z,-}} (z, z^\prime)
\end{pmatrix}. \notag
\end{align}
It follows that
\begin{equation}
\tr \left( \hD e^{-s^2 {\hD}^2} 
((u, z), (u^\prime,z^\prime)) \right) =
i \: \frac{1}{\sqrt{4\pi s^2}} \: \frac{u-u^\prime}{2s^2} \: 
e^{-\frac{(u - u^\prime)^2}{4s^2}} \: \tr_s \left( 
e^{-s^2 D_Z^2}(z, z^\prime) \right).
\end{equation}
Hence
\begin{align}
\eta_{\langle k \rangle}(M) & = \frac{2}{\sqrt{\pi}} \: \int_0^\infty
\int_Z i \:  \frac{1}{\sqrt{4\pi s^2}} \: \frac{k}{2s^2} \: 
e^{-\frac{k^2}{4s^2}} \: \tr_s \left( 
e^{-s^2 D_Z^2}(\phi^k(z), z) \right) d\vol(z) \: ds\\
& = \frac{2}{\sqrt{\pi}} \: \int_0^\infty
i \:  \frac{1}{\sqrt{4\pi s^2}} \: \frac{k}{2s^2} \: 
e^{-\frac{k^2}{4s^2}} \: \Tr_s \left( 
\phi^k e^{-s^2 D_Z^2} \right) \: ds \notag \\
& = \frac{2}{\sqrt{\pi}} \: \int_0^\infty
i \:  \frac{1}{\sqrt{4\pi s^2}} \: \frac{k}{2s^2} \: 
e^{-\frac{k^2}{4s^2}} \: \Tr_s \left( 
\phi^k \Big|_{Ker(D_Z)} \right) \: ds \notag \\
& = \frac{i}{k\pi}  \: \Tr_s \left( 
\phi^k \Big|_{Ker(D_Z)} \right) \: \int_0^\infty
e^{-\frac{k^2}{4s^2}} \: d \left( - \frac{k^2}{4s^2} \right) \notag \\
& = \frac{i}{k\pi}  \: \Tr_s \left( 
\phi^k \Big|_{Ker(D_Z)} \right). \notag 
\end{align}
The proposition follows. \qed

\section{Proof of Proposition \ref{toshow8}} \label{proof8}

Any irreducible unitary representation $\rho$ of 
$F \widetilde{\times}_\alpha \Z$ arises
as follows \cite[Section 10]{Mackey (1978)}. 
First, $\Z$ acts on the dual space $\hF$.  A periodic point
of period $j$ corresponds to a representation $\mu : F \rightarrow 
U(N)$ and a matrix $U \in U(N)$ such that $\mu(\alpha^j (f)) =
U \mu(f) U^{-1}$. (The matrix $U$ is determined up to multiplication by
a unit complex number.) 
Consider the representation $\nu : F \widetilde{\times}_\alpha j \Z
\rightarrow U(N)$
given by
\begin{equation}
\nu (f,jr) = \mu(f) \: U^r.
\end{equation}
Then $\rho$ comes from inducing $\nu$ from 
$F \widetilde{\times}_\alpha j \Z$ to
$F \widetilde{\times}_\alpha \Z$.  
The character of $\rho$ is
\begin{equation} \label{character}
\chi_\rho(f,k) = 
\begin{cases}
0 &\text{if $j \nmid k$,} \\
\Tr \left( \left[ \mu(f) + \mu(\alpha^{-1}(f)) + \ldots +
\mu(\alpha^{-(j-1)}(f)) \right] U^r \right) &\text{if $k=jr$}. \\   
\end{cases}
\end{equation}

Let ${\cal T}_\rho(M)$ be the analytic torsion of $M$, computed using
the representation $\rho$.
From Fourier analysis,
\begin{equation} \label{us1}
{\cal T}_\rho(M) = \sum_{f,k} \chi_\rho(f,k) \: {\cal T}_{\langle f,k \rangle}
(M).
\end{equation}
Let $M^\prime$ be the mapping torus of
$\phi^j$. Then $M^\prime$ is a $j$-fold cover of $M$. By
\cite[(VI), p. 27]{Fried (1987)}, 
\begin{equation} \label{us2}
{\cal T}_\rho(M) = {\cal T}_{\nu}(M^\prime).
\end{equation} 

Let $E_\mu$ be the flat $\C^N$-bundle on $Z$ coming from the representation
$\mu$. Then $\phi^j$ acts on $Z$ and preserves $E_\mu$. Let
$L_\mu(r)$ be the Lefschetz number of $\phi^{jr}$ acting on
$(Z, E_\mu)$. Put
\begin{equation}
\zeta_\nu(z) = \exp \left( \sum_{r>0} \frac{z^r}{r} \: L_\mu(r) \right).
\end{equation}
It follows from
\cite[Section 3]{Milnor (1968)} that
\begin{equation} \label{into}
{\cal T}_{\nu} (M^\prime) = \ln |\zeta_\nu(1)|^2.
\end{equation}

Take a cellular decomposition of $Z$. Let $\hZ$ have the lifted cellular
structure and let $C^*(\hZ)$ denote the cellular cochains on $\hZ$.
We let $F$ act on $C^*(\hZ)$ on the right by
\begin{equation}
\omega \cdot f = R_{f^{-1}}^* \omega.
\end{equation}
Then
\begin{equation} \label{equiv}
\hphi^* (\omega \cdot f) = (\hphi^* \omega) \cdot \alpha^{} (f).
\end{equation}
We can identify $C^*(Z; E_\mu)$ with $C^*(\hZ) \otimes_F \C^N$, with
the relation $\omega \cdot f \otimes_F v = \omega \otimes_F \mu(f) v$.
Then $\phi^{jr}$ acts on $C^*(\hZ) \otimes_F \C^N$ by
\begin{equation}
\phi^{jr} \left( \omega \otimes_F v \right) = \left( \hphi^{jr} \right)^*
\omega \otimes_F U^{r} V.
\end{equation}

Letting $\Tr_s$ denote the supertrace on $C^*(Z; E_\mu)$, we want to
compute 
\begin{equation} \label{tocompute}
L_\mu(r) = \Tr_s \left( \phi^{jr} \right).
\end{equation}
For the moment we
concentrate on $C^p(Z; E_\mu)$.
Let $\{e_i\}$ be a basis of $C^p(\hZ)$ consisting of dual $p$-cells.
The set of such dual $p$-cells has a free $F$-action.
Write the action of $\hphi^*$ on $C^p(\hZ)$ as
\begin{equation}
\hphi^* ( e_i) = \sum_l 
\hphi^*_{e_i \rightarrow e_l} e_l.
\end{equation}
From (\ref{equiv}),
\begin{equation}
\hphi^*_{e_i \rightarrow e_l} =
\hphi^*_{e_i f \rightarrow e_l \alpha^{} (f)}.
\end{equation}

Let $\{\overline{e}_i\}$ be a set of representatives for the $F$-orbits of
the dual $p$-cells.
Then as a vector space,
\begin{equation}
C^p(\hZ) \otimes_F \C^N = \bigoplus_i \overline{e}_i \otimes \C^N
\end{equation}
and so
\begin{align}
\phi^{jr} \left( \overline{e}_i \otimes v \right) & =
\sum_{l,f} \left(\hphi^{jr}\right)^*_{ \overline{e}_i 
\rightarrow \overline{e}_l f} \overline{e}_l f \otimes_F U^{r}v \\
& =
\sum_{l,f} \left(\hphi^{jr}\right)^*_{ \overline{e}_i 
\rightarrow \overline{e}_l f} \overline{e}_l \otimes \mu(f) U^{r}v. \notag
\end{align}
Hence
\begin{equation}
\Tr \left( \phi^{jr} \right) = \sum_{i,f}  
\left(\hphi^{jr}\right)^*_{ \overline{e}_i \rightarrow \overline{e}_i f}
\Tr \left(\mu(f) U^{r} \right).
\end{equation}
As the choice of the representatives $\{ \overline{e}_i \}$ is arbitrary, 
we can also write
\begin{equation} \label{use1}
\Tr \left( \phi^{jr} \right) = \frac{1}{|F|} \sum_{i,f}  
\left(\hphi^{jr}\right)^*_{ e_i \rightarrow e_i f}
\Tr \left(\mu(f) U^{r} \right).
\end{equation}

We have
\begin{align} \label{use2}
\Tr \left( \phi^{jr} \right) & = \frac{1}{|F|} \sum_{i,f}  
\left(\hphi^{jr}\right)^*_{ e_i \rightarrow e_i f}
\Tr \left(\mu(f) U^{r} \right) \\
& = \frac{1}{|F|} \sum_{i,l,f}  
\left(\hphi \right)^*_{ e_i \rightarrow e_l}
\left(\hphi^{jr-1}\right)^*_{ e_l \rightarrow e_i f}
\Tr \left(\mu(f) U^{r} \right) \notag \\
& = \frac{1}{|F|} \sum_{i,l,f}  
\left(\hphi^{jr-1}\right)^*_{ e_l \rightarrow e_i f}
\left(\hphi \right)^*_{ e_i \rightarrow e_l}
\Tr \left(\mu(f) U^{r} \right) \notag \\
& = \frac{1}{|F|} \sum_{i,l,f}  
\left(\hphi^{jr-1}\right)^*_{ e_l \rightarrow e_i f}
\left(\hphi \right)^*_{ e_i f \rightarrow e_l \alpha^{} (f)}
\Tr \left(\mu(f) U^{r} \right) \notag \\
& = \frac{1}{|F|} \sum_{l,f}  
\left(\hphi^{jr}\right)^*_{ e_l \rightarrow e_l \alpha^{} (f)}
\Tr \left(\mu(f) U^{r} \right) \notag \\
& = \frac{1}{|F|} \sum_{l,f}  
\left(\hphi^{jr}\right)^*_{ e_l \rightarrow e_l f}
\Tr \left(\mu(\alpha^{-1}(f)) U^{r} \right). \notag 
\end{align}
Then from (\ref{use1}) and (\ref{use2}),
\begin{equation} \label{thenfrom}
\Tr \left( \phi^{jr} \right) = \frac{1}{j|F|} \sum_{i,f}  
\left(\hphi^{jr}\right)^*_{ e_i \rightarrow e_i f}
\Tr \left( \left[\mu(f) + \mu(\alpha^{-1}(f)) + \ldots + 
\mu(\alpha^{-(j-1)}(f)) \right] U^{r} \right).
\end{equation}

Put 
\begin{equation} \label{np}
n_{p,jr}(f) = \frac{1}{|F|} \sum_{i,f}  
\left(\hphi^{jr}\right)^*_{ e_i \rightarrow e_i f}.
\end{equation}
From (\ref{character}), (\ref{tocompute}) and (\ref{thenfrom}),
\begin{equation}
L_\mu (r) = \frac{1}{j} \sum_f \chi_\rho(f,jr) \sum_{p = 0}^n (-1)^p \: 
n_{p,jr}(f).
\end{equation}

Put
\begin{equation}
i_{p,jr}(f) = \sum_{f^\prime \sim_{jr} f} \sum_{i}  
\left(\hphi^{jr}\right)^*_{ \overline{e}_i \rightarrow 
\overline{e}_i f}.
\end{equation}
The Nielsen fixed-point index 
$I_{jr}(f) \in \Z$ of the transformation $\phi^{jr}$ 
is defined by \cite[Section 1]{Geoghegan-Nicas (1994)}
\begin{equation}
I_{jr}(f) = \sum_{p = 0}^n (-1)^p \: i_{p,jr}(f).
\end{equation}

Put $s_{jr}(f) = | \{ \gamma \in F : 
\gamma f \alpha^{jr}(\gamma^{-1}) = f\}|$.
We have
\begin{align} \label{ip}
i_{p,jr}(f) & = \frac{1}{s_{jr}(f)} \sum_{i, \gamma} 
\left(\hphi^{jr}\right)^*_{ \overline{e}_i \rightarrow 
\overline{e}_i \gamma f \alpha^{jr}(\gamma^{-1})} \\   
& = \frac{1}{s_{jr}(f)} \sum_{i, \gamma} 
\left(\hphi^{jr}\right)^*_{ \overline{e}_i \gamma \rightarrow 
\overline{e}_i \gamma f} \notag \\
& = \frac{1}{s_{jr}(f)} \sum_{i} 
\left(\hphi^{jr}\right)^*_{ e_i \rightarrow e_i f}. \notag
\end{align}
Then from (\ref{np}) and (\ref{ip}),
\begin{equation}
n_{p,jr}(f) = \frac{s_{jr}(f)}{|F|} \: i_{p,jr}(f) = 
\frac{i_{p,jr}(f)}{|[f]_{jr}|}.
\end{equation}
Hence the Lefschetz number is given in terms of the Nielsen index by
\begin{equation} \label{LN}
L_\mu(r) = \frac{1}{j} \sum_f 
\frac{\chi_\rho(f,jr) \: I_{jr}(f)}{|[f]_{jr}|}.
\end{equation} 
Substituting (\ref{LN}) into (\ref{into}) and using
(\ref{us1})-(\ref{us2}) gives (\ref{eq8}).
As $\Gamma$ is a type-I discrete group \cite[p. 61]{Mackey (1978)},
knowing equation (\ref{eq8}) for all $\rho \in \widehat{\Gamma}$ 
determines $\{{\cal T}_{\langle f,k \rangle}(M)\}_{(f,k) \in \Gamma}$. \qed

\section{Examples} \label{Examples}
\noindent
Proposition \ref{toshow1} : It follows from Proposition \ref{finite}.1 that
$b_{p, \langle g \rangle}(M)$ can be nonzero if $\Gamma$ is finite.
For example, $b_{0, \langle g \rangle}(M) = 
\frac{|\langle g \rangle|}{|\Gamma|}$. \\ \\
\noindent
Proposition \ref{toshow2} : It follows from Proposition \ref{finite}.2 that
${\cal T}_{\langle g \rangle}(M)$ is nonzero in some examples in which $M$ is
a lens space.\\ \\
\noindent
Proposition \ref{toshow3} : It follows from Proposition \ref{finite}.3 that
$\eta_{\langle g \rangle}(M)$ is nonzero in some examples in which $M$ is
a lens space, both for the tangential signature operator and the Dirac
operator.\\ \\
\noindent
Proposition \ref{toshow4} : An even-dimensional closed hyperbolic manifold
$M$ satisfies the hypotheses of the proposition for all $p$. \\ \\
\noindent
Proposition \ref{toshow5} : Let $N_1$ be a closed even-dimensional spin 
manifold whose fundamental group is virtually nilpotent or Gromov-hyperbolic,
with $\widehat{A}(N_1) \ne 0$.
Let $N_2$ be a lens space which is spin and whose Dirac operator has a
nonzero $\rho$-invariant. Put $M = N_1 \times N_2$. Shrink $N_2$ so that
$M$ has positive scalar curvature. Then $M$ satisfies the
hypotheses of the proposition and
$\eta_{\langle e, g \rangle}(M) = \widehat{A}(N_1) \:
\eta_{\langle g \rangle}(N_2)$ is nonzero for appropriate $g$. \\ \\
\noindent
Proposition \ref{toshow6} : There are many nontrivial examples.\\ \\
\noindent
Proposition \ref{toshow7} : Nontrivial examples come from 
closed even-dimensional
oriented manifolds $Z$ with a finite-order orientation-preserving 
diffeomorphism $\phi$ such
that $\phi$ has nonzero Lefschetz or Atiyah-Bott numbers. \\ \\
\noindent
Proposition \ref{toshow8} : Any example of Proposition \ref{toshow7} gives
an example of Proposition \ref{toshow8} by taking $F$ to be the trivial
group. There are also many examples with $F$ nontrivial.  

\end{document}